\documentclass[11pt]{amsart}
\usepackage{amssymb,amsmath,amsthm}
\usepackage{a4}
\usepackage{epsf}
\usepackage{epsfig}
\begin{document}

\newcommand{\NP}{$\mathcal{N}\mathcal{P}$}
\newcommand{\newsymb}{{\mathcal P}}
\newcommand{\Pd}{{\mathcal P}^d}
\newcommand{\R}{{\mathbb R}}
\newcommand{\N}{{\mathbb N}}
\newcommand{\QQ}{{\mathbb Q}}
\newcommand{\ZZ}{{\mathbb Z}}
\newcommand{\STAB}{\mathrm{STAB}}

\newcommand{\inter}{\mathrm{int}}
\newcommand{\conv}{\mathrm{conv}}
\newcommand{\aff}{\mathrm{aff}}
\newcommand{\lin}{\mathrm{lin}}
\newcommand{\dist}{\mathrm{dist}}
\newcommand{\trans}{\intercal}
\newcommand{\diam}{\mathrm{diam}}
\newcommand{\pp}{\mathfrak{p}}
\newcommand{\pf}{\mathfrak{f}}
\newcommand{\pg}{\mathfrak{g}}
\newcommand{\PP}{\mathfrak{P}}
\newcommand{\pl}{\mathfrak{l}}
\newcommand{\pv}{\mathfrak{v}}
\newcommand{\cl}{\mathrm{cl}}
\newcommand{\bx}{\overline{x}}

\def\ip(#1,#2){#1\cdot#2}

\newtheorem{theorem}{Theorem}[section]
\newtheorem{corollary}[theorem]{Corollary}
\newtheorem{lemma}[theorem]{Lemma}
\newtheorem{remark}[theorem]{Remark}
\newtheorem{definition}[theorem]{Definition}  
\newtheorem{proposition}[theorem]{Proposition}  
\newtheorem{claim}[theorem]{Claim}
\numberwithin{equation}{section}

\title[Representation of  polyhedra by polynomial inequalities]
      {On the representation of polyhedra by polynomial inequalities} 


\author{Martin Gr\"otschel}
\address{Martin Gr\"otschel, Konrad-Zuse-Zentrum f\"ur
  Informationstechnik (ZIB), 
         Taku\-str.~7, D-14195 Berlin-Dahlem, Germany}
\email{groetschel@zib.de}

\author{Martin Henk}
\address{Martin Henk, Technische Universit\"at Wien, Abteilung f\"ur Analysis,
  Wiedner Hauptstr. 8-10/1142, A-1040 Wien, Austria} 
\curraddr{Technische Universit\"at Berlin, Institut f\"ur Mathematik,
  Sekr. MA 6-2, Stra{\ss}e des 17. Juni 136 
D-10623 Berlin, Germany}
\email{henk@tuwien.ac.at}

\begin{abstract}
A beautiful  result of Br\"ocker and Scheiderer on the stability index
of basic closed semi-algebraic sets implies, as a very special case, that
every $d$-dimensional polyhedron admits a representation as the set of
solutions of at most $d(d+1)/2$ polynomial inequalities. Even
in this polyhedral case, however, no constructive proof
is known, even if the 
quadratic upper bound is replaced by any bound depending only on the
dimension. 

Here we give, for simple polytopes, 
an explicit construction of polynomials describing such a
polytope. The number of used polynomials is  exponential in the
dimension, but  in the 2- and 3-dimensional case we get the
expected number $d(d+1)/2$. 
\end{abstract}

\maketitle

\section{Introduction}
By a surprising  and striking result of Br\"ocker and Scheiderer (see
\cite{Scheiderer:stability}, \cite{broecker:semialgebraic}, 
 \cite{bochnakcosteroy:real_algebraic_geometry}  and \cite{Mahe:broecker_scheiderer}) every basic closed semi-algebraic set  of the form      
\begin{equation*}
      \mathcal{S}=\left\{ x\in \R^d : \pf_1(x)\geq 0,\dots, \pf_l(x)\geq 0\right\}, 
\end{equation*}
where $\pf_i\in \R[x_1,\dots,x_d]$, $1\leq i\leq l$,  are
polynomials, can be represented by at most $d(d+1)/2$ polynomials,
i.e., there exist polynomials 
$\pp_1,\dots,\pp_{d(d+1)/2}\in\R[x_1,\dots,x_d]$ such that  
 \begin{equation*} 
     \mathcal{S}=
     \left\{ x\in \R^d : \pp_1(x)\geq 0,\dots, \pp_{d(d+1)/2}(x)\geq 0\right\}.
\end{equation*}
Moreover, in the case of basic open semi-algebraic sets, i.e., $\geq$ is
replaced by strict inequality, one can even bound the maximal number
of needed polynomials by the dimension $d$ instead of $d(d+1)/2$.  

No explicit constructions, however, of such  
systems of polynomials are known, even in the very special case of
$d$-dimensional convex polyhedra  and even if the 
quadratic upper bound is replaced by any bound depending only on the
dimension. In \cite{broecker:semialgebraic}, Example
2.10, or in \cite{andradas_broecker_ruiz:constructible}, Example 4.7,
a description of a regular convex $m$-gon in the plane 
by two polynomials is given. This result was generalized to arbitrary
convex polygons and three polynomial inequalities by vom Hofe
\cite{vom_Hofe:disser}.   
Bernig \cite{Bernig:diplom} proved that, for $d=2$,  every  convex polygon
can  even be represented by two polynomial inequalities.  
The main purpose of this note is to give  some basic properties of
polynomials describing polyhedra as well as an explicit  construction of
(exponentially many) polynomials describing simple $d$-polytopes of
any dimension $d$.

In order to state the result we fix some notation. 
A polyhedron $P\subset \R^d$ is the intersection of 
finitely many closed halfspaces, i.e., we can represent it as 
\begin{equation}
  P=\left\{x\in\R^d : \ip(a^i,x) \leq b_i, \, 1\leq i\leq
    m\right\},
\label{eq:poly}
\end{equation} 
for some $a^i\in\R^d$, $b_i\in\R$. Here $\ip(a,x)$ denotes the
standard 
inner product on $\R^d$. If $P$ is bounded then it is called
a polytope. In general we are only working with
$d$-dimensional polyhedra  $P\subset\R^d$, and for
short, we denote these polyhedra as $d$-polyhedra. 
A $d$-polyhedron $P$ 
is called simple, if every
$k$-dimensional face, $0\leq k\leq d-1$,  belongs to exactly
$d-k$ facets  of $P$.  In the case of polytopes, this is equivalent to
the statements that every vertex lies in precisely $d$ facets, or,
every vertex figure is a simplex (cf.~\cite{ziegler:polytopes}, pp.~65).
Since for unbounded  polyhedra the above definition of simple polyhedra  is
not invariant with respect to projective transformation, we call a
polyhedron $P$ a {\em strongly simple} polyhedron,  
if it is projectively equivalent to a
simple polytope,   
i.e., for a ``$P$-permissible'' projective transformation $f:\R^d\to\R^d$ the 
closure of $f(P)$, denoted by $\cl(f(P))$, is a simple polytope.    
For more information about polyhedra, polytopes and their 
combinatorial structure we refer to  the books \cite{McMShe:upper} and 
\cite{ziegler:polytopes}.
 
For polynomials $\pp_i\in \R[x_1,\dots,x_d]$, $1\leq i\leq l$, we
denote by  
\begin{equation*}
 \newsymb(\pp_1,\dots,\pp_l) := \left\{ x\in \R^d : \pp_1(x)\geq
   0,\dots, \pp_{l}(x)\geq 0\right\} 
\end{equation*}
the associated closed semi-algebraic set and we define 
\begin{definition} A {\em $\newsymb$-representation} of a $d$-polyhedron
  $P\subset\R^d$ consists of $l$ polynomials
  $\pp_1,\dots,\pp_l\in \R[x_1,\dots,x_d]$ such that  
\begin{equation*}
 P=\newsymb\left(\pp_1,\dots,\pp_l\right).
\end{equation*}
\end{definition}
For $d$-polytopes there are two other well known and
important representations (see, e.g., \cite{GritzmannKlee:handbook_crc},
\cite{ziegler:polytopes}):  The representation of $P$ by  $m$ vectors
$a^i\in\R^d$ and scalars $b_i$ as given in
\eqref{eq:poly} is called {\em $\mathcal{H}$-representation}
of $P$. Of course, any $\mathcal{H}$-representation may be
regarded as a special $\newsymb$-representation of $P$ with linear
polynomials (linear forms). As a
dual counterpart we have the 
 {\em $\mathcal{V}$-representation} 
of a $d$-polytope $P$ consisting of $n$ points $v^i\in\R^d$ such that 
$P$ is the convex hull of these points, i.e., $P=\conv\{v^1,\dots,v^n\}$. 

Both, $\mathcal{V}$- and $\mathcal{H}$-representations, are quite
powerful and useful representations of  polytopes. They admit  
the   computation of  the complete combinatorial
structure (face-lattice) of the polytope
(cf.~\cite{GritzmannKlee:handbook_crc}, \cite{Seidel:computational})
and linear programming problems can
be solved in polynomial time with respect to the input sizes of these representations.  Many interesting combinatorial
  optimization problems, however, cannot be effectively solved via
  these representations because  the size of both
  representations is exponential in the ``natural'' input size of the
  combinatorial problem instances. This holds, e.g., for the polytopes associated
  with the traveling salesman problem or the max cut problem, see
  \cite{GLS:ellipsoid} for details.   

On the other hand, the result of Br\"ocker and Scheiderer
  tells us that there always 
  exists  a $\newsymb$-representation by polynomially  many (with
  respect to the dimension) polynomials, and
  therefore, a 
  representation or ``good'' approximation of a polytope by
  few polynomial inequalities could lead to  a new
  approach to ``hard'' combinatorial optimization  problems via nonlinear
  programming tools. In the last section we discuss this connection in
  more detail as well as the 
  possible outcomes of such an approach.
 
  For a different approach
  of approximating  
  ``discrete problems'' by semi-alge\-braic sets see
  \cite{Barvinok_vershik:semi_approximation}, and for related problems
  involving polynomials and optimization see, e.g., 
  \cite{barvinok:norm_poly}, 
  \cite{brieden_gritzmann:poly_obj}, 
  \cite{lasserre:moments} and the references within.

Unfortunately, at the moment our knowledge about polynomials
representing
or approximating polytopes is rather limited. For
arbitrary polytopes we even do not know 
how to convert -- via an algorithm --  a
$\mathcal{H}$-representation into a $\mathcal{P}$-representation where the
number of polynomials depends only on the dimension.  For simple
polytopes we  have the following result.

\begin{theorem} Let $P\subset \R^d$ be a $d$-dimensional simple
  polytope given by a $\mathcal{H}$-representation. 
  Then $\mu(d)\leq d^d$ polynomials
  $\pp_i\in\R[x_1,\dots,x_d]$ can be constructed such that 
\begin{equation*}
   P =\newsymb\left(\pp_1,\dots,\pp_{\mu(d)}\right).
\end{equation*}
In particular, we can take $\mu(2)=3$ and $\mu(3)=6$.
\label{thm:main}
\end{theorem}
Since every $2$-dimensional polygon is
simple, Theorem \ref{thm:main} includes the result of 
 vom Hofe \cite{vom_Hofe:disser}. 

It is not hard to see  that if
a polyhedron is given  
as the set of solutions of polynomial inequalities then the sum of the total
degrees of these polynomials is at least the number of facets (cf.~Proposition \ref{prop:basic} i)). In
fact, the total degrees of the polynomials used in Theorem
\ref{thm:main} depend on the number of  $k$-faces,
$k=0,\dots,d-1$, of the polytope as well as on some metric
properties of the polytope.  For upper bounds on the degrees in the
general semi-algebraic setting we refer to \cite{buresimahe:polynomial_bound}.

All the polynomials that we use in Theorem \ref{thm:main} are of a
rather special structure, namely, if $P=\{x\in\R^d : 
\ip(a^i,x)\leq b_i, 1\leq i\leq m\}$ then they can be expressed as 
\begin{equation*}
  \sum_{\alpha\in\N^m,\,\alpha\geq 0}c_\alpha \prod_{i=1}^m
  \left[b_i-\ip(a^i,x)\right]^{\alpha_i},  
\label{eq:handelman_type}
\end{equation*}
where $c_\alpha$ are certain non-negative numbers and, of course, only
finitely many of them are positive.  One possible explanation 
for this special type is given by a result of Handelman
\cite{handelman:polynomials} which
says that every polynomial that is strictly positive on  a polytope $P$ is
of that type. His proof is non-constructive, for a more explicit
version see \cite{powersreznick:polya}. 

It seems to be an interesting question to ask for the minimal number
of polynomials needed to describe a
$d$-polyhedron or $d$-polytope. Therefore we define 

\begin{definition} For a $d$-polyhedron $P\subset \R^d$, let $m_\newsymb(P)$
  be the minimal number of polynomials needed in a
  $\newsymb$-representation of $P$ and let 
\begin{equation*}
\begin{split}
\overline{m}_\newsymb(d)&:=\max\left\{m_\newsymb(P): P\subset \R^d \text{ is a
      $d$-polyhedron}\right\}, \\
  m_\newsymb(d)&:=\max\left\{m_\newsymb(P): P\subset \R^d \text{ is a
      $d$-polytope}\right\}.
\end{split}
\end{equation*}
\end{definition}
We set $m_\newsymb(\R^d)=0$, and for a polyhedron $\tilde{P}\subset
\R^d$ with $\dim(\tilde{P})<d$ we mean by $m_\newsymb(\tilde{P})$ the
minimal number of polynomials in $\dim(\tilde{P})$-variables, which
are needed in order to describe an embedding of $\tilde{P}$ in
$\R^{\dim(\tilde{P})}$. 

Observe that $m_\newsymb(P)$ is invariant with respect to regular affine
transformations of $P$. Moreover, 
it is easy to see that $m_\newsymb(P)\geq d$ for every $d$-polytope
$P$  
(cf.~Corollary \ref{cor:basic} i)). 
 And together with the result  of Br\"ocker and Scheiderer we obtain 
\begin{equation*}
      d \leq m_\newsymb(d)\leq  \overline{m}_\newsymb(d)\leq d(d+1)/2.
\end{equation*}
In Proposition \ref{prop:polytope_polyhedra} we
show  $\overline{m}_\newsymb(d)\leq m_\newsymb(d)+1$.   
Probably, the truth is  $\overline{m}_\newsymb(d)= m_\newsymb(d)$.  
 
 There
are some trivial examples of polytopes for which $d$ polynomials 
are sufficient. 
For instance, the cube $C^d=\{ x\in \R^d : -1\leq x_i \leq
1\}$ can be written as $C^d=\{ x\in \R^d : -(x_i)^2+1\geq
0\}$. Another example is an arbitrary $d$-simplex $T^d$. To see
this, we may assume w.l.o.g.~that $T^d=\{x\in\R^d : x_i\geq 0,\,
x_1+\cdots +x_d\leq 1\}$. Then it is easy to check that 
$$
  T^d=\left\{x\in\R^d : x_i\left(1-x_i-\cdots -x_d\right)\geq
    0,\,1\leq i\leq d\right\}. 
$$ 
Actually, the given representations of a cube and a simplex are
special cases of a general construction of polynomial inequalities for
prisms and pyramids (cf.~Proposition \ref{prop:prism_pyramid}), which
in particular imply that every $3$-dimensional prism or pyramid can be
described by $3$ polynomials (cf.~Corollary \ref{cor:3_pyramid}). 
However, we are not aware of a representation
of a regular cross\-polytope 
$C^d_\star=\{x\in\R^d : |x_1|+\cdots + |x_d|\leq 1 \}$ by $d$
polynomials or of any constructive good upper bound on
$m_\newsymb(C^d_\star)$. 

In this context it seems to be worth mentioning that a classical
result of Minkowski \cite{minkowski:volumen} implies that a polytope
can be approximated 
``arbitrarily closely'' by only one  polynomial
inequality, where --  of course -- the degree of this polynomial is
``arbitrarily large''.
In Section 2 we will construct such a
polynomial, which will be  used in the
scope of the proof of Theorem \ref{thm:main}.  
Furthermore, in Section 2 we will state some simple and basic properties
of polynomials describing a polyhedron. In Section 3 we
give the construction of the polynomials used in  Theorem
\ref{thm:main} and the proof of this theorem. A generalization of the
theorem to strongly simple polyhedra is given in Section 4 (cf.~Corollary
\ref{cor:polyhedra}). Finally, in section 5 we discuss possible
outcomes of research on  $\newsymb$-representation of  polyhedra
associated with hard
combinatorial optimization problems.


\section{Polynomials and Polyhedra}

Let 
\begin{equation}
  P=\left\{x\in\R^d : \ip(a^i,x) \leq b_i,\, 1\leq i\leq m  \right\}
\label{eq:polytope_structure}
\end{equation}
be a $d$-dimensional polyhedron. 
We always assume that  the representation
\eqref{eq:polytope_structure} is irredundant, i.e., $P\cap\{x\in\R^d :
\ip(a^i,x) =
b_i\}$ is a facet of $P$, $1\leq i\leq m$. In particular, 
we have that $b_i=h(a^i)$, where $h(\cdot)$ denotes the support
function of $P$, i.e., 
$$
   h(u)=\sup\left\{\ip(u,x) : x\in P\right\}.
$$ 
For a non-negative linear combination of vectors $\sum_i\rho_i\,u^i$,
$\rho_i\in\R_{\geq 0}$, we have 
$h(\sum_i\rho_i\,u^i)\leq\sum_i\rho_i h(u^i)$. 
The next proposition collects some simple properties of polynomials
describing polyhedra. 

\begin{proposition} Let $P=\{x\in\R^d : \ip(a^i,x) \leq b_i,\,
  1\leq i\leq m\}$ be a $d$-polyhedron and let
$\pp_1,\dots,\pp_l\in\R[x_1,\dots,x_d]$ such that $P=\newsymb(\pp_1,\dots,\pp_l)$.

\begin{enumerate}
  \item[{\rm i)}] Each facet defining linear polynomial $b_i-\ip(a^i,x)$,
     $i\in\{1,\dots,m\}$, is a factor of one of the  $\pp_j$.
   \item[{\rm ii)}] Let $F$ be a $k$-dimensional face of $P$. Then
    there exist $d-k$ polynomials $\pp_{j_1},\dots,\pp_{j_{d-k}}$,
    say, such that these polynomials vanish on the affine hull of $F$,
    i.e., $$\aff F\subset \left\{ x\in \R^d :  
        \pp_{j_1}(x)=\cdots =\pp_{j_{d-k}}(x)=0 \right\}.$$ 
\end{enumerate}
\label{prop:basic}
\end{proposition}
\begin{proof} {\rm i)}\quad 
Let $F_i=P\cap \{x\in \R^d : \ip(a^i,x) =
  b_i\}$. By assumption, $F_i$ is a facet of $P$.   First we note that
  for each $y\in F_i$ one of the polynomials $\pp_j$ has to
  vanish. Otherwise,  if $\pp_j(y)>0$ for all $1\leq j\leq l$ we
  get by the continuity of polynomials that we can move $y$ in
  any direction without leaving $P$, which contradicts the property
  that $y$ belongs to the boundary. Thus we know that the
  polynomial  
$$
      \mathfrak f(x) = \prod_{j=1}^l \pp_j(x)  
$$ 
vanishes on $F_i$. Hence it vanishes on every segment joining two
points of $F_i$ and therefore, it has to be $0$ on
$\aff(F_i)=\{x\in\R^d : b_i-\ip(a^i,x)=0\}$. Thus  
 $b_i-\ip(a^i,x)$ is a factor of $\mathfrak{f}(x)$. 
Furthermore, since
$b_i-\ip(a^i,x)$ is irreducible, it has to be a factor of one of
the $\pp_j$ (see, e.g., \cite{CoxLittleOShea:ideals_etc}, pp.~148). 

\noindent{\rm ii)}\quad We use induction with respect to the dimension
$k$ of the 
face $F$ and we
  start with $k=d-1$. In this case the statement follows immediately
  from i). So let $k<d-1$ and let $G$ be a
$(k+1)$-face containing $F$. By induction we can assume that 
$\aff(G)\subset \{x\in\R^d : \pp_1(x)=\pp_2(x)=\cdots =\pp_{d-k-1}(x)=0\}$ and
so 
\begin{equation*}
  G=\aff(G)\cap P=\{x\in\aff(G) : \pp_{d-k}(x)\geq
  0,\dots,\pp_l(x)\geq 0\}.
\end{equation*}
With respect to the $(k+1)$-dimensional polytope $G$ in the space
$\aff(G)$ the face $F$ is a facet 
and so we can conclude that one of the 
polynomials $\pp_{d-k},\dots,\pp_l$ vanishes on $\aff(F)$. 
\end{proof}
As an immediate consequence of Proposition
\ref{prop:basic} ii) we note
\begin{corollary} \hfill
\begin{enumerate}
\item[{\rm i)}]  Let $F$ be an non-empty face of a $d$-polyhedron
  $P\subset\R^d$. Then 
\begin{equation*}
   m_\newsymb(P)\geq m_\newsymb(F)+d-\dim(F).
\end{equation*}
In particular, $m_\newsymb(P)\geq d$ for $d$-polytopes.
\item[{\rm ii)}] $m_\newsymb(d+1)\geq m_\newsymb(d)+1$ and 
   $\overline{m}_\newsymb(d+1)\geq \overline{m}_\newsymb(d)+1$. 
\end{enumerate}
\label{cor:basic}
\end{corollary}
\begin{proof} Let $P=\newsymb(\pp_1,\dots,\pp_l)$. By Proposition
  \ref{prop:basic} ii) we may assume that $\aff(F)\subset\{x\in\R^d
  : \pp_1(x)=\cdots=\pp_{d-\dim(F)}(x)=0\}$. Thus 
\begin{equation*}
  F=\aff(F)\cap P=\{x\in\aff(F) : \pp_{d-\dim(F)+1}(x)\geq
  0,\dots,\pp_l(x)\geq 0\} 
\end{equation*}
and so $m_\newsymb(F)\leq m_\newsymb(P) -(d-\dim(F))$. If $P$ is a
polytope then we may choose for $F$ a vertex and get
$m_\newsymb(P)\geq d$.

  For ii) let
  $Q$ be a $d$-polytope  with $m_\newsymb(Q)=m_\newsymb(d)$. Now
  we take  any $(d+1)$-polytope that has $Q$ as a facet and then we
  can conclude from i) that  $m_\newsymb(P)\geq m_\newsymb(d)+1$. Of
  course, the polyhedral case can be treated analogously.   
\end{proof}
The next statement gives some  information on  $m_\newsymb(P)$ for
$d$-prisms and $d$-pyramids. A $d$-polytope $P$ is called a {\em
  $d$-pyramid} 
({\em $d$-prism}) with basis $Q$, where $Q$ is a $(d-1)$-polytope,
if there exists a $v\in\R^d$ such that $P=\conv\{Q,v\}$ 
($P=Q+\conv\{0,v\}=\{q+\lambda\,v : q\in Q,\,0\leq\lambda\leq 1\}$).  

\begin{proposition} Let $P$ be a $d$-prism or a $d$-pyramid with basis
  $Q$. Then 
\begin{equation*}
         m_\newsymb(P)=m_\newsymb(Q)+1.
\end{equation*}
\label{prop:prism_pyramid}
\end{proposition}
\begin{proof} 
Since in both cases the $(d-1)$-polytope $Q$ is a facet of $P$  we get
from Corollary \ref{cor:basic} i)
\begin{equation}
             m_\newsymb(P)\geq m_\newsymb(Q)+1.
\label{eq:lower_bd}
\end{equation} 
In order to show the reverse inequality we start with a
$d$-dimensional pyramid
$P=\conv\{Q,v\}$ and w.l.o.g.~we assume  that 
\begin{equation}
  {\rm i)}\,\,
   Q\subset \left\{x\in\R^d : x_d=0 \text{ and }\sum_{i=1}^{d-1}
   (x_i)^2\leq 1\right\},\quad {\rm ii)}\,\,
   v=(0,\dots,0,1)^\trans.
\label{eq:assump_prop}
\end{equation}
Let $l=m_\newsymb(Q)$ and $\pp_1,\dots,\pp_l\in\R[x_1,\dots,x_{d-1}]$ such 
 that $Q=\newsymb(\pp_1,\allowbreak\dots,\pp_l)$. 
Furthermore we denote by $p$ the maximum of the total 
 degrees of the polynomials $\pp_j$, $1\leq j\leq l$, and let 
\begin{equation}
\begin{split}
 \tilde{\pp}_j(x)&=(1-x_d)^p\cdot\pp_j\left(\frac{x_1}{1-x_d},\dots,\frac{x_{d-1}}{1-x_d}\right),\quad 1\leq j\leq l,\\
        \pp(x)&=x_d\left(1-x_d-x_d\left((x_1)^2+\cdots +(x_{d-1})^2\right)\right).
\end{split}
\label{eq:last_poly}
\end{equation}
Observe that $\tilde{\pp}_j(x)\in\R[x_1,\dots,x_d]$, $1\leq j\leq l$. Next we claim that 
\begin{equation}
  P=\newsymb(\tilde{\pp}_1,\dots,\tilde{\pp}_l,\,\pp).     
 \label{eq:prep}
\end{equation}
To see this we first note that for $0\leq \lambda \leq1$    
\begin{equation}
  P\cap\{x\in\R^d : x_d=\lambda\} = (1-\lambda)\,Q+\lambda\, (0,\dots,0,1)^\trans.
\label{eq:pyramid}
\end{equation}
A simple calculation shows that, for $0\leq \lambda<1$, 
\begin{equation}
\begin{split}
 P &\cap\{x\in\R^d :  x_d =\lambda\} \\ &=\left\{(x_1,\dots,x_{d-1},\lambda)^\trans\in\R^d  :
    \tilde{\pp}_j(x_1,\dots,x_{d-1},\lambda)\geq 0,\, 
    1\leq j\leq l\right\}.
\end{split}
\label{eq:section}
\end{equation}
Next we observe that, for $x\in P$, we have $0\leq x_d\leq 1$. By
\eqref{eq:pyramid} and  assumption \eqref{eq:assump_prop} i) we
conclude that 
$$
        (x_1)^2+\cdots (x_{d-1})^2 \leq (1-x_d)^2 \text{ for all }
        x\in P.
$$
Hence $\pp(x)\geq 0$ for  $x\in P$ and together with
\eqref{eq:section} we get 
$P\setminus\{(0,\dots,0,1)^\trans\}\subset
\newsymb(\tilde{\pp}_1,\dots,\tilde{\pp}_l,\,\pp)$   
and consequently $P\subset
\newsymb(\tilde{\pp}_1,\dots,\tilde{\pp}_l,\,\pp)$.

For the reverse inclusion we notice that $\pp(x)\geq 0$ 
implies $0\leq x_d\leq 1$ and with  \eqref{eq:section} we obtain 
\begin{equation*}
   \newsymb(\tilde{\pp}_1,\dots,\tilde{\pp}_l,\,\pp)\setminus\{x\in\R^d : x_d=1\} 
    \subset P. 
\end{equation*}
Since for $x_d=1$ the inequality $\pp(x)\geq 0$ becomes 
$(x_1)^2+\cdots +(x_{d-1})^2 \leq 0$ we conclude that 
$$
   \newsymb(\tilde{\pp}_1,\dots,\tilde{\pp}_l,\,\pp)\cap\{x\in\R^d :
   x_d=1\}=\{(0,\dots,0,1)^\trans\}.  
$$
Hence we also have
$\newsymb(\tilde{\pp}_1,\dots,\tilde{\pp}_l,\,\pp)\subset P$. Thus  
\eqref{eq:prep} is shown and so we have $m_\newsymb(P)\leq
m_\newsymb(Q)+1$. Together with $\eqref{eq:lower_bd}$ the statement of
the proposition is verified for pyramids.

If $P=Q+\conv\{0,v\}$ is a $d$-prism over the basis $Q$ and if we
assume again that $v=(0,\dots,0,1)^\trans$, $Q\subset \{x\in\R^d:
x_d=0\}$ and $Q=\newsymb(\pp_1,\dots,\pp_{m_\newsymb(Q)})$ then it is easy to check that 
\begin{equation*}
\begin{split}
  P =\left\{x\in\R^d: \pp_j(x_1,\dots,x_{d-1})\geq 0, \, 1\leq j\leq
m_\newsymb(Q),\, x_d(1-x_d)\geq 0\right\}.
\end{split}
\end{equation*}  
\end{proof}
Since every $2$-dimensional polygon can be described by two
polynomials (cf.~\cite{Bernig:diplom}), Proposition \ref{prop:prism_pyramid} gives  
\begin{corollary} Let $P$ be a $3$-dimensional prism or pyramid. Then 
\begin{equation*}
             m_\newsymb(P)=3.
\end{equation*}
\label{cor:3_pyramid}
\end{corollary}
\vspace{-0.5cm}
Next we study the relation between $m_\newsymb(d)$ and
$\overline{m}_\newsymb(d)$. Obviously, we have $m_\newsymb(d)\leq
\overline{m}_\newsymb(d)$. In order to bound
$\overline{m}_\newsymb(d)$ in terms of $m_\newsymb(d)$, we apply
a standard technique from Discrete Geometry, which ``makes an unbounded
pointed polyhedron bounded'', namely projective transformations. 
\begin{proposition} Let $d\geq 2$. Then 
$$
m_\newsymb(d)\leq \overline{m}_\newsymb(d)\leq m_\newsymb(d)+1.
$$
\label{prop:polytope_polyhedra}
\end{proposition}
\vspace{-0.5cm}
\begin{proof} In order to prove the
  upper bound on $\overline{m}_\newsymb(d)$ let $P$ be a $d$-polyhedron
  such that $\overline{m}_\newsymb(d)=m_\newsymb(P)$.  Let $G$ be a non-empty
  face of minimal dimension of $P$,  and we assume that
$0\in G$. Suppose that $\dim(G)>0$.  Then 
the intersection of  $P$ with the orthogonal complement of $\lin(G)$,  
the linear hull of $G$,  is a lower dimensional polyhedron $Q$, say. 
Since $P=Q+\lin(G)=\{q+g : q\in Q,\, g\in\lin(G)\}$ any 
$\newsymb$-representation of  $Q$ can easily be converted  to a
$\newsymb$-representation of $P$ with the same number of polynomials.
 With the help of Corollary \ref{cor:basic} ii) we get the contradiction 
\begin{equation*} 
   \overline{m}_\newsymb(d)=m_\newsymb(P)= m_\newsymb(Q)
   \leq \overline{m}_\newsymb(\dim(Q)) < \overline{m}_\newsymb(d). 
\end{equation*}
Therefore, we can assume that  the origin is 
a vertex of $P$.  
Thus, we can find a vector $c\in\R^d$
with $\ip(c,x) >0$ for all $x\in P\setminus\{0\}$. Let
$f:\R^d\to\R^d$  be the projective map
$$
          f(x)=\frac{x}{\ip(c,x) +1}.
$$
Then we can describe $f(P)$ by a set of inequalities of the form 
$f(P)=\{x\in\R^d : Ax\leq b,\, \ip(c,x) <1\}$, for a certain matrix
$A\in\R^{m\times d}$ and a vector $b\in\R^m$. 
The inequality $\ip(c,x) <1$ corresponds to the points at infinity. 

We conclude that  the set (the closure of $f(P))$ 
\begin{equation}
 \cl(f(P))=\{x\in\R^d : Ax\leq b,\, \ip(c,x) \leq 1\}
\label{eq:closure_polyhedra}
\end{equation}
 is a $d$-dimensional polytope. Hence  we get can find a
 $\newsymb$-representation of $\cl(f(P))$ by polynomials $\pp_i$,
 $i\in I$, say, with $\#I\leq
m_\newsymb(d)$. So we may write 
\begin{equation*}
f(P)=\left\{x\in\R^d : \pp_i(x)\geq 0,\, i\in I,\,\, \ip(c,x) <
  1\right\}.
\label{eq:polyhedra_2}
\end{equation*}
Thus  
$$ 
    P=\left\{x\in\R^d\setminus\{x\in \R^d : \ip(c,x)=-1\} : \pp_i( f(x) )\geq 0,\, i\in I,\,\,\ip(c,f(x)) <1\right\}.
$$
Since $\ip(c,f(x))<1$ is equivalent to $\ip(c,x) +1 > 0$ we may
multiply all rational functions $\pp_i( f(x) )$ by suitable powers of
$\ip(c,x) +1$ and obtain some polynomials $\tilde{\pp}_i(x)$,
say, such that 
$$
  P=\left\{x\in\R^d : \tilde{\pp}_i( x )\geq 0,\, i\in I,\,\,
     \ip(c,x) +1 >0\right\}.
$$
Since $\ip(c,x)\geq 0$ for all $x\in P$ we may replace the last
inequality in this representation by $\ip(c,x) \geq 0$ and since $\#I\leq
m_\newsymb(d)$ the
proposition is shown.
\end{proof}

In the next lemma  a strictly convex polynomial  $\pp$ is constructed such
that the convex body $K=\{x\in \R^d : \pp(x)\leq 1\}$  is not 
``too far away'' from $P$. Here  the distance between
convex bodies $K_1,K_2$ will be measured by  the Hausdorff distance
$\dist(K_1,K_2)$, i.e.,
$$
  \dist(K_1,K_2):=\max\left\{ \max_{x\in K_1}\,\min_{y\in K_2}\Vert
    x-y\Vert,\, \max_{x\in K_2}\,\min_{y\in K_1}\Vert
    x-y\Vert\right\},
$$
where $\Vert\cdot\Vert$ denotes the Euclidean norm. Furthermore, for a
bounded set $S\subset \R^d$, the diameter is denoted by $\diam(S)$, i.e., 
$$
      \diam(S):=\max\{\Vert x-y\Vert : x,y\in S\}.
$$
In order to construct this strictly convex polynomial we follow an 
approach of Hammer \cite{hammer:minkowski}, but since we need a
slightly different approximation we give the short
proof. For similar results see  \cite{firey:approximating} and
\cite{weil:schachteln}. 

\begin{lemma} Let $P=\{x\in\R^d : \ip(a^i,x) \leq b_i,\,\, 1\leq i\leq
  m\}$ be a $d$-dimensional polytope. For $1\leq i\leq m$ let 
$$ 
  \pv_i(x):=
   \frac{2\ip(a^i,x) - h(a^i)+h(-a^i)}
        {h(a^i)+h(-a^i)}.
$$
Let $\epsilon>0$, $p >
\ln(m)/(2\ln(1+\frac{2\epsilon}{(d+1)\diam(P)}))$, 
$$
  \pp_\epsilon(x):=
 \sum_{i=1}^m \frac{1}{m}\,\left[\pv_i(x)\right]^{2\,p} \quad\text{ and }\quad
  K_\epsilon:=\{x\in\R^d : \pp_\epsilon(x)\leq 1\}.
$$
Then we have 
$P\subset K_\epsilon$ and $\dist(P,K_\epsilon)\leq\epsilon$.
\label{lem:approximating} 
\end{lemma} 
\begin{proof} Since $|\pv_i(x)|\leq 1$ for all $x\in P$ we certainly
  have $P\subset K_\epsilon$. W.l.o.g. let the origin be the center of
  gravity of $P$ and let $\lambda=\epsilon/\diam(P)$. First we check
  that $K_\epsilon\subset P_\lambda=\{x\in\R^d : \ip(a^i,x) \leq
  (1+\lambda)\,b_i,\,1\leq i\leq m\}$. 
  Let $y\notin P_\lambda$. Then we may assume
  $\ip(a^1,y) > (1+\lambda)h(a^1)$ which implies 
$$
   \pv_1(y)> 1+2\lambda \frac{h(a^1)}{h(a^1)+h(-a^1)} 
   \geq 1 +\frac{2\lambda}{d+1},
$$
where the last inequality follows from the choice of the origin as the
center of gravity (cf.~\cite{BonFen:con}, p.52). By the lower bound on $p$
we conclude $(1/m)\,\pv_1(y)^{2p}>1$ and thus $\pp_\epsilon(y)>1$,
which shows $y\notin K_\epsilon$. Finally we observe that
$\dist(P,P_\lambda)\leq \lambda\, \diam(P)=\epsilon$. 
\end{proof}

\section{Proof of Theorem \ref{thm:main}}

In the following let 
\begin{equation*}
  P=\left\{x\in\R^d : \ip(a^i,x) \leq b_i,\, 1\leq i\leq m  \right\}
\end{equation*}
be a convex  $d$-dimensional {\em simple} polytope with $m$ facets. We
further assume that we know all $k$-faces of the polytope as well as 
the facets containing a given face. This information can be
obtained from the $\mathcal{H}$-representation above by several
(exponential time and space) algorithms (cf.~\cite{Seidel:computational}). 

We remark
that every $k$-face of $P$ is contained in exactly $d-k$ facets.  
The set of all $k$-dimensional faces is denoted by 
$\mathcal{F}_k$, $0\leq k\leq d-1$. 
For a $k$-face $F$ of $P$, let $[F]_1,\dots,[F]_{d-k}$,
$[F]_1\leq\cdots \leq [F]_{d-k}$, be 
all indices of vectors $a^i$ such $\ip(a^i,x) = b_i$,
for all $x\in\aff(F)$. In other words, these are the ordered indices
of all facets containing $F$.  

Next, for a $k$-face $F$ and a positive integral  vector $w\in\N^{d-k}$,
we define 
\begin{equation}
     a(F,w):=\sum_{j=1}^{d-k} w_j\, a^{[F]_j}.  
\label{eq:support_vectors}
\end{equation}
Observe that $a(F,w)$ is a support vector of $F$, i.e., 
\begin{equation*}
      F=P\cap\left\{x\in\R^d : \ip({a(F,w)}, x) = h(a(F,w))\right\}.
\end{equation*}
Moreover, since $F$ is contained in all the facets corresponding to
the vectors $a^{[F]_j}$ we note that 
\begin{equation}
 h(a(F,w))=\sum_{j=1}^{d-k} w_j\, h(a^{[F]_j}).
\label{eq:supp_func}
\end{equation}

With    $w\in\N^{d-k}$ and the set $\mathcal{F}_k$ of all $k$-faces we associate the polynomial 
\begin{equation}
  \pp_{k,w}(x):=\prod_{F\in \mathcal{F}_k} [h(a(F,w))-\ip({a(F,w)},x)].
\label{eq:polynomials_can}
\end{equation}
So, for a fixed $w$ the polynomial $\pp_{k,w}(x)$ is the product of
all those supporting hyperplanes of all $k$-faces of $P$  which can be
written as in \eqref{eq:support_vectors}.  Since we are only interested in finitely many different support
vectors of the type $a(F,w)$ at a given face $F$ we define certain
sets of integral vectors:
\begin{equation*}
 \mathcal{W}_{d-1}:=\{(1)\},\quad
 \mathcal{W}_{d-k}:=\{(2^{l_1},\dots,2^{l_k})^\trans  : 0\leq l_i\leq
 k-2\},\quad 2\leq k\leq d.
\end{equation*}
In particular we have $\mathcal{W}_{d-2}=\{(1,1)\}$ and 
\begin{equation}
       \#\mathcal{W}_{k}=(d-k-1)^{d-k}.
\label{eq:card_setM}
\end{equation} 
The meaning of these sets $\mathcal{W}_k$ will be explained  in the
next lemma. 
\begin{lemma} Let $P$ be a simple $d$-dimensional polytope. Let 
  $k\geq 1$, $F,G\in\mathcal{F}_k$ with  $F\cap G\ne\emptyset$ and let
  $w\in\mathcal{W}_k$, $y\in\R^d$ such that 
\begin{equation}
      h(a(F,w))-\ip({a(F,w)},y) < 0 \quad\text{and}\quad 
      h(a(G,w))-\ip({a(G,w)},y) \leq 0.
\label{eq:assump_lemma}
\end{equation}
Then there exists  an
$\tilde{w}\in\mathcal{W}_{\dim(F\cap G)}$ such that  
\begin{equation*}
h(a(F\cap G,\tilde{w}))-\ip({a(F\cap G,\tilde{w})}, y) < 0.
\end{equation*} 
\label{lem:set_M}
\end{lemma}
\begin{proof} In view of \eqref{eq:supp_func} 
  we conclude from 
  \eqref{eq:assump_lemma} that 
\begin{equation}
\begin{split}
   \sum_{j=1}^{d-k} w_j
       \left(h(  a^{[F]_j} )-\ip(a^{[F]_j},
         y)\right)&<0, \\
  \quad 
   \sum_{j=1}^{d-k} w_j
       \left(h(  a^{[G]_j} )-\ip(a^{[G]_j},
         y)\right)
  &\leq 0.
\end{split}
\label{eq:conclude_fom_assump}
\end{equation} 
Let $I_1=\{[{F}]_1,\dots,[{F}]_{d-k}\}$,
$I_2=\{[{G}]_j : h(a^{[{G}]_j})-\ip(a^{[{G}]_j},
         y) \leq 0\}$,  
$I_3=\{[{G}]_j : h(a^{[{G}]_j})-\ip(a^{[{G}]_j},
y) > 0\}$ and let $I=I_1\cup I_2\cup I_3$.

Further, we need  a map $\tau$ that gives, for a $k$-face
$F$ and a number $q\in\{[F]_1,\dots,[F]_{d-k}\}$, the position
of $q$ with respect to  the ordered list $[F]_1,\dots,\allowbreak[F]_{d-k}$. 
In particular we have  $\tau(F,[F]_j)=j$.  Now we merge the two normal vectors
$a(F,w), a(G,w)$ in the following way: Let   
\begin{equation*}
\begin{split}
  a &= \sum_{j\in I_1\cap I_2}
      2\,\max\{w_{\tau(F,j)},w_{\tau(G,j)} \}a^j+ 
      \sum_{j\in I_1\cap I_3} a^j \\ 
    &+ \sum_{j\in (I_2\cup I_3)\setminus I_1} w_{\tau(G,j)}a^j + 
       \sum_{j\in I_1\setminus (I_2\cup I_3)} w_{\tau(F,j)}a^j.
\end{split}
\end{equation*}
On account of \eqref{eq:conclude_fom_assump} we have 
\begin{equation*}
   h(a)- \ip(a,y)<0.
\end{equation*}
Since the polytope is simple, the assumption 
${F}\cap{G}\ne\emptyset$,  implies that the vectors $\{a^j : j\in I\}$
are the vectors of all facets containing the face $F\cap G$ and  
by construction we may write $a=\sum_{j\in I} \tilde{w}_{i_j}a^j$
for some numbers $\tilde{w}_{i_j}=2^{l_{i_j}}$ with $0\leq l_{i_j}\leq
d-k-1$.  

Thus we have $a=a(F\cap G,\tilde{w})$ for a certain vector $\tilde{w}
\in\mathcal{W}_{\dim(F\cap G)}$. 
\end{proof}

We note that from the proof of Lemma \ref{lem:set_M} it follows that
it suffices to define the set $\mathcal{W}_{d-3}$ as 
\begin{equation}
      \mathcal{W}_{d-3}=\left\{(1,1,2)^\trans,\, (1,2,1)^\trans,\,
        (2,1,1)^\trans \right\}.
\label{eq:set_M3}
\end{equation}

Lemma  \ref{lem:set_M} says that if two linear factors of a
polynomial $\pp_{k,w}$, $k\in\{1,\dots,d-1\}$, $w\in\mathcal{W}_{k}$  
are non-positive and at least one is negative, 
than there exists a linear factor of a polynomial of the type
$\pp_{\tilde{k},\tilde{w}}$, $\tilde{k}< k$,
  $\tilde{w}\in\mathcal{W}_{\tilde{k}}$, which has to be negative,
  too. Therefore, with these sets $\mathcal{W}_k$ we associate the
following sets of polynomials: 
\begin{equation}
    \PP_{k}=\left\{ \pp_{k,w}(x) : w \in\mathcal{W}_k\right\}, 
   \quad k=0,\dots,d-1. 
\label{eq:polyset}
\end{equation} 
$\PP_{d-1}$, $\PP_{d-2}$ consist of only one polynomial,
$\#\PP_{d-3}=3$ and for
$0\leq k\leq d-3$ we have (cf.~\eqref{eq:card_setM})
\begin{equation}
             \# \PP_k =(d-k-1)^{d-k}.
\label{eq:poly_card}
\end{equation}

 We need one more polynomial. To this  end we set for two
vectors $a,b\in\R^d\setminus\{0\}$ 
\begin{equation*}
  U(a,b)=\left\{x\in\R^d : \ip(a,x) \geq h(a) \text{ and } \ip(b,
  x) \geq h(b)\right\}.
\end{equation*}
$U(a,b)$ is a closed set and so we can define 
\begin{equation*}
  \epsilon(a,b)=\min\left\{ \Vert x-y \Vert : x\in P,\, y\in U(a,b) \right\}.
\end{equation*} 
Since 
$P\subset \left\{x\in\R^d : \ip(a,x) \leq h(a),\,\ip(b,x) \leq h(b)\right\}$ and both planes $\{x\in\R^d : \ip(a,
x)=h(a)\}$, $\{x\in\R^d : \ip(b,x)=h(b)\}$ are supporting
hyperplanes we have 
\begin{equation}
  \epsilon(a,b)=0 \Longleftrightarrow 
   \left\{x\in\R^d : \ip(a,x) = h(a) \text{ and } \ip(b,
  x) = h(b)\right\}\cap \mathcal{F}_0\ne\emptyset.
\label{eq:dist_vert} 
\end{equation}
For $0\leq k\leq d-1$ we set 
\begin{equation}
  \epsilon_{k}  = \min\left\{ \epsilon\left(a(F,w),a(\tilde{F},w)\right)>0 :
    F\ne\tilde{F}\in\mathcal{F}_k,\, w\in \mathcal{W}_k\right\}.
\label{eq:epsilon_k}
\end{equation}
We note that for different vertices $v,\tilde{v}\in\mathcal{F}_0$, $w\in\mathcal{W}_0$, we
always have (cf.~\eqref{eq:dist_vert}) 
\begin{equation}
      \epsilon\left(a(v,w),a(\tilde{v},w)\right)>0.
\label{eq:diff_vert}
\end{equation}
Finally  let $\overline{\epsilon}$ satisfying  
\begin{equation}
  0 < \overline{\epsilon} <\min\left\{ \epsilon_k : 0\leq k\leq d-1 \right\}.
\label{eq:epsilon}
\end{equation}
 With respect to $\overline{\epsilon}$ let
 $\pp_{\overline{\epsilon}}(x)$ be the polynomial according to Lemma
 \ref{lem:approximating} and let  
\begin{equation*}
\begin{split}
       \mathcal{V}_1(P) &= \left\{ x \in \R^d :  \PP_k(x)\geq 0, \,1\leq
         k\leq d-1,\,\, \pp_{\overline{\epsilon}}(x)\leq 1 \right\}, \\
       \mathcal{V}_2(P) &= \left\{ x \in \R^d : 
                          \PP_{0}(x)\geq 0  \right\}.
\end{split}
\end{equation*}
Here for a set of polynomials $\PP$, say, $\PP(x)\geq 0$ means
$\pp(x)\geq 0$ for all $\pp(x)\in \PP$. 

Before giving the last piece of the proof of Theorem \ref{thm:main} we
remark that in order to find a  number $\overline{\epsilon}$ satisfying
\eqref{eq:epsilon} we have to calculate several
distances $\epsilon(a,b)$. In general, $\epsilon(a,b)$ can be
calculated  (or sufficiently well approximated) by several LP-based methods
 (cf.~\cite{MatouSharirWelzl:lp}). 
 In particular, depending on the input size of the polytope and
 the vectors $a,b$ one 
 can give a lower bound on this distance if it is positive. Thus for
 a given polytope we can calculate such an $\overline{\epsilon}$ and
 hence   the polynomial $\pp_{\overline{\epsilon}}(x)$. 

\begin{proof}[Proof of Theorem \ref{thm:main}] On account of
  \eqref{eq:poly_card}  and \eqref{eq:set_M3} 
  the  theorem will follow from the identity  
\begin{equation}
  P= \mathcal{V}_1(P)\cap \mathcal{V}_2(P).
\label{eq:main_equation}
\end{equation} 
Obviously, by the definition of all these polynomials via support
vectors and by Lemma \ref{lem:approximating} we know that $P$ is
contained in the set on the right hand side.  
In order to prove the reverse inclusion  we first claim 
\begin{claim} Let $y\notin P$, but $y\in
  \mathcal{V}_1(P)$.
  Then there exists a vertex $v$ of $P$ and an $w\in\mathcal{W}_0$   
  such that $h(a(v,w))-\ip({a(v,w)},y) < 0$.
\label{claim:claim}
\end{claim}
Since $y\notin P$ at least one of the
inequalities  $\ip(a^i, x)\leq b_i$,
  $i\in\{1,\dots,m\}$, is violated and so we may define $k$ as the smallest
  dimension such that there exists a face $F\in\mathcal{F}_k$ 
  and  an $w\in \mathcal{W}_{k}$
  with $h(a(F,w))-\ip({a(F,w)},y) < 0$. Suppose $k>0$. 
  Since $y\in \mathcal{V}_1(P)$ we have $p_{k,w}(y)\geq 0$ and so 
  there must exist another $k$-face $G$ with 
  $h(a(G,w))- \ip({a(G,w)},y)\leq 0$. Hence we have 
\begin{equation*}
       y\in U(a(F,w),a(G,w)).
\end{equation*}
If $\epsilon(a(F,w),a(G,w))>0$ then we get from the definition of
$\overline{\epsilon}$  and the approximating
property of the polynomial $\pp_{\overline{\epsilon}}(x)$ (cf. Lemma
\ref{lem:approximating}) that $\pp_{\overline{\epsilon}}(y)>1$. Thus
we can assume that $\epsilon(a(F,w),a(G,w))=0$ and  
 from \eqref{eq:dist_vert} we get $F\cap G\ne\emptyset$. 
Therefore we may apply Lemma \ref{lem:set_M} and we get a
contradiction of 
the minimality of $k$. This shows Claim \ref{claim:claim}.

Now let $y\notin P$. We want to show that $y$ is not contained in the
set $\mathcal{V}_1(P)\cap\mathcal{V}_2(P)$. Suppose that $y\in
\mathcal{V}_1(P)\cap\mathcal{V}_2(P)$. By Claim
\ref{claim:claim} we may assume that there exists a vertex
$v\in\mathcal{F}_0$ and an $w\in\mathcal{W}_0$ such that 
 $h(a(v,w))-\ip({a(v,w)},y) < 0$. However,
 $y\in \mathcal{V}_2(P)$ implies $\pp_{0,w}(y)\geq 0$ and thus
 there exists another vertex $\tilde{v}$ with
 $h(a(\tilde{v},w))-\ip({a(\tilde{v},w)},y) \leq
 0$. Therefore we have 
\begin{equation*} 
            y\in W\left(a(v,w),a(\tilde{v},w)\right).
\end{equation*}
Next we observe that $\epsilon(a(v,w),a(\tilde{v},w))>0$
(cf.~\eqref{eq:diff_vert}) and by the definition of
$\overline{\epsilon}$ we conclude $\pp_{\overline{\epsilon}}(y)>1$,
which gives the contradiction $y\notin \mathcal{V}_1(P)$.

\end{proof}

\section{Remarks}
First we want to generalize Theorem \ref{thm:main} to the class of
strongly simple polyhedra. As in the case of polytopes, a
$\mathcal{H}$-representation of a polyhedron $P$ is a description of
$P$ by linear inequalities of the form \eqref{eq:poly}

\begin{corollary} Let $P\subset \R^d$ be a $d$-dimensional strongly simple  
                  polyhedron given by a $\mathcal{H}$-representation. 
Then $\mu(d)\leq d^d$ polynomials $\pp_i\in\R[x_1,\dots,x_d]$ 
can be constructed such that 
\begin{equation*}
   P =\left\{ x\in \R^d :  \pp_i(x)\geq 0,\,\,\, 1\leq i\leq \mu(d)\right\}. 
\end{equation*}
\label{cor:polyhedra}
\end{corollary}
\begin{proof} The proof is just a combination of the proofs 
  of Proposition \ref{prop:polytope_polyhedra} and Theorem
  \ref{thm:main}.  As in the proof 
  of Proposition \ref{prop:polytope_polyhedra} we first note that we
  can assume that $P$ has a vertex. Next we apply a projective
  transformation $f(x)=x/(\ip(c,x)+1)$ such that $\cl(f(P))$
  becomes a polytope. By the definition of strongly simple polyhedra,  
  $\cl(f(P))$  is a $d$-dimensional simple polytope. 
  Hence, from Theorem \ref{thm:main}, we get a
  representation   of the type 
$$
f(P)=\left\{x\in\R^d : \pp_i(x)\geq 0,\, i\in I,\,\, \ip(c,x) <
  1\right\}
$$
with certain polynomials $\pp_i(x)$, $i\in I$, $\#I\leq \mu(d)$. Now we
can proceed as in the proof of Proposition
\ref{prop:polytope_polyhedra} in order to get a $\newsymb$-representation of $P$
with $\mu(d)+1$ polynomials.  A closer look on the number  $\mu(d)$
shows that $\mu(d)<d^d$ for $d\geq 2$ (cf.~\eqref{eq:poly_card})  and so the assertion is proved. 
\end{proof}

We remark that the  proofs of Theorem \ref{thm:main} and of Corollary \ref{cor:polyhedra} can be adapted such
that for  arbitrary polyhedra a 
representation  by  polynomials is obtained where the number of
polynomials depends  exponentially   on the dimension and the maximal
degree of a vertex of the polytope. In other words, degeneracy in the
sense of linear programming leads to additional difficulties.  
Since, however, the main problem is 
to find a representation of a polytope by a few polynomials we omit
a proof of this statement.

In the $2$-dimensional case the meaning of the polynomials $\pp_1(x)$,
$\pp_0(x)$ and $\pp_{\overline{\epsilon}}(x)$ can easily be
illustrated. Suppose the polygon is the $7$-gon depicted in these  
pictures.
\begin{figure}[htb]
\begin{center}
\hbox{\input{2poly_first.pstex_t}\hfill\input{2poly_second.pstex_t}}
\end{center}
\begin{center}
 \mbox{$\{x\in\R^d : \pp_1(x)\geq
     0\}$}\hspace{4cm}\mbox{$\{x\in\R^d : \pp_0(x)\geq 0\}$}
\end{center}
\caption{}
\label{fig:one}
\end{figure}    
Then the shaded region  on the left hand side of Figure \ref{fig:one}
shows all points in the plane that 
satisfy the inequality $\pp_1(x)\geq 0$, whereas the shaded regions on
the right hand side correspond to the points $\pp_0(x)\geq 0$. 
If we intersect the shaded regions of both pictures  we get the points
satisfying both inequalities (see Figure \ref{fig:two}).

\begin{figure}[htb]
\begin{center}
\input{2poly_remain.pstex_t}
\end{center}
\begin{center}
\mbox{$\{x\in\R^d : \pp_1(x)\geq 0\, \text{and}\,\pp_0(x)\geq 0\}$}
\end{center}
\caption{}
\label{fig:two}
\end{figure}

We see that all the shaded points that do not belong to the  polygon are
``far away'' from the polygon and thus we can cut them off with the
inequality $\pp_{\overline{\epsilon}}(x)\leq 1$.

\begin{figure}[htb]
\vspace{-1cm}
\epsfig{file=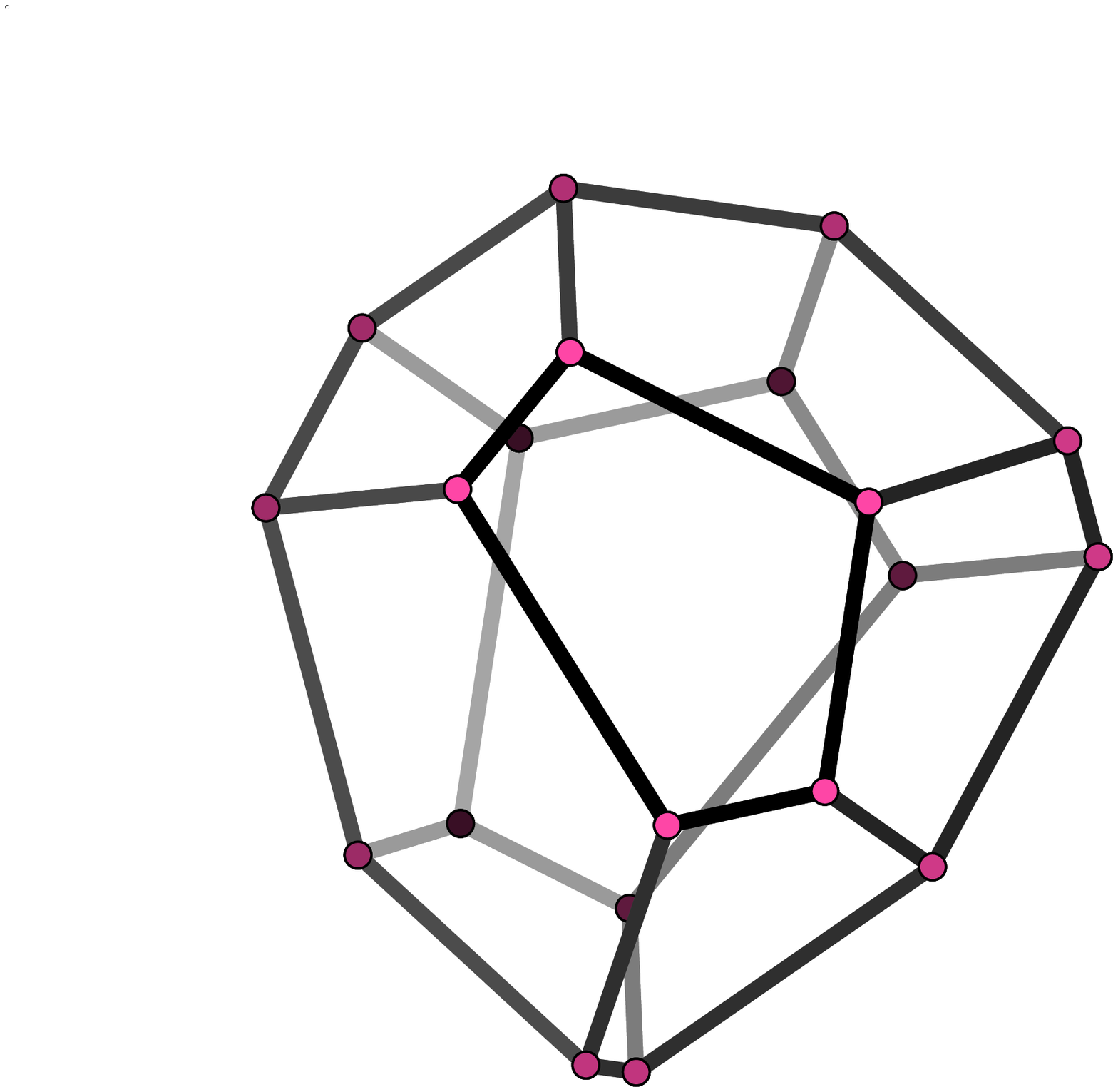, width=6truecm}

\vspace{-0.5cm}
\centering{{\sc Figure 3.} ({\small produced using polymake \cite{polymake} and  
          javaview \cite{javaview}}) }
\label{fig:three}
\end{figure}

Now  let $P=\{x \in \R^3 : \ip(a^i,x) \leq b_i,\,1\leq i\leq m\}$
be a simple $3$-dimensional convex polytope. With the notation from
section 3  we get the following polynomials (cf.~\eqref{eq:set_M3})
\begin{equation*}
 \begin{split}
     \pp_{2}(x) & = \prod_{i=1}^m \left[b_i-\ip(a^i,x) \right], \\
     \pp_{1}(x) & = \prod_{F\in\mathcal{F}_1} 
            \left[(b_{[F]_1}+b_{[F]_2})-\ip({(a^{[F]_1}+a^{[F]_2})},x)
            \right], \\
     \pp_{0,(1,1,2)}(x) & = \prod_{v\in\mathcal{F}_0} 
            \left[(b_{[v]_1}+b_{[v]_2}+2\,b_{[v]_3})-
                  \ip({(a^{[v]_1}+a^{[v]_2}+2\,a^{[v]_3})},x)\right], \\
     \pp_{0,(1,2,1)}(x) & = \prod_{v\in\mathcal{F}_0} 
            \left[(b_{[v]_1}+2\,b_{[v]_2}+b_{[v]_3})-
                  \ip({(a^{[v]_1}+2\,a^{[v]_2}+a^{[v]_3})},x)\right],
                \\
      \pp_{0,(2,1,1)}(x) & = \prod_{v\in\mathcal{F}_0} 
            \left[(2\,b_{[v]_1}+b_{[v]_2}+b_{[v]_3})-
                  \ip({(2\,a^{[v]_1}+a^{[v]_2}+a^{[v]_3})},x)\right], \\
    \pp_{\overline\epsilon}(x) & =
                     \sum_{i=1}^m \frac{1}{m}\,
   \left[\frac{2\ip(a^i,x) -
       b_i+h(-a^i)}{b_i+h(-a^i)}\right]^{2\,p}, 
 \end{split}
\end{equation*}  
where $\overline{\epsilon}$  has to be chosen such that
\eqref{eq:epsilon} is satisfied, and $p$ is given by Lemma
\ref{lem:approximating}. Let us consider a ``real''
3-dimensional polytope $P=\{x\in\R^d : Ax\leq b\}$, with 
{\small
\begin{equation*}     
       A=\left(\begin{array}{rrr}
          0 & 3 & 2 \\
          0 & -3 & 2 \\
          2 &  0 & 3 \\
          2 &  0 & -3 \\
          3 & 2 & 0 \\
         -3 & 2 & 0 \\
          0 & -3 & -2 \\
          0 & 3 & -2 \\
          -2 & 0 & -3 \\
          -2 & 0 & 3 \\
          -3 & -2 & 0 \\
          3 & -2  & 0 
         \end{array}\right),\quad\quad 
       b=\left(\begin{array}{r} 5\\6\\ 5\\4\\5\\5\\6\\5\\6\\5\\4\\6
          \end{array}\right).
\end{equation*}
} 
$P$ is a simple polytope with  12 facets, all of them  pentagons, 30 edges, 20 vertices, and it
may be described  as a ``skew'' dodecahedron (see Figure \ref{fig:three}). 

With respect to  the facets we get
the polynomial
{\tiny
\begin{equation*}
 \begin{split}
   \pp_2(x)=&(5-3x_2-2x_3)(6+3x_2-2x_3)(5-2x_1-3x_3)(4-2x_1+3x_3)\\
            &(5-3x_1-2x_2)(5+3x_1-2x_2)(6+3x_2+2x_3)(5-3x_2+2x_3)\\
            &(6+2x_1+3x_3)(5+2x_1-3x_3)(4+3x_1+2x_2)(6-3x_1+2x_2).
 \end{split}
\end{equation*} 
}
For the 30 edges we  obtain
{\tiny
\begin{equation*}
\begin{split}
  \pp_1(x)=
 &(10-2x_1-3x_2-5x_3)(10+  2x_1-3x_2-5x_3)(10+3x_1+-5x_2-2x_3)\\
 &(10-6x_2)(10-3x_1-5x_2-2x_3)(11-2x_1+  3x_2-5x_3)\\
 &(12-3x_1+5x_2-2x_3)(12+6x_2)(10+  3x_1+  5x_2-2x_3)\\
 &(11+  2x_1+  3x_2-5x_3)(10-6x_3)(11-5x_1+  2x_2-3x_3)\\
 &(10-5x_1+-2x_2-3x_3)(9-5x_1-2x_2+  3x_3)(10-5x_1+  2x_2+  3x_3)\\
 &(10-2x_1+  3x_2+  5x_3)(10+  6x_3)(9-2x_1-3x_2+  5x_3)\\
 &(10-3x_1-5x_2+  2x_3)(11-6x_1)(10+  5x_1-2x_2-3x_3)\\
 &(9+  6x_1)(11+  5x_1-2x_2+  3x_3)(10+  3x_1-5x_2+  2x_3)\\
 &(12-3x_1+  5x_2+  2x_3)(12+  2x_1+  3x_2+  5x_3)(10+  3x_1+  5x_2+
 2x_3)\\
 &(11+  2x_1+-3x_2+  5x_3)(10+  5x_1+  2x_2+  3x_3)(9+  5x_1+  2x_2-3x_3).
\end{split}
\end{equation*} 
}
With respect to the 20 vertices we get 3 polynomials depending on the
weights $w\in\mathcal{W}_0=\{(1,1,2),(1,2,1),(2,1,1)\}$.
{\tiny
\begin{equation*}
\begin{split}
  \pp_{0,(1,1,2)}(x)=
 &(20-8x_1-7x_2-5x_3)(20+  2x_1-3x_2-11x_3)(20-3x_1+
 -11x_2+  2x_3)\\
&(20+  3x_1-11x_2+  2x_3)(20+  7x_1-5x_2-8
 x_3)(21+  2x_1+  3x_2-11x_3)\\
&(23-8x_1+  7x_2-5x_3)
 (20+  6x_1+  10x_2)(24-6x_1+  10x_2)\\
&(19+  8x_1+  7x_2+
 -5x_3)(22-11x_1+  2x_2-3x_3)(19-5x_1-8x_2+  7x_3)\\
&(21+-11x_1+  2x_2+  3x_3)(22+  2x_1+  3x_2+  11x_3)(22-8x_1+  7x_2
 +  5x_3)\\
&(21+  2x_1-3x_2+  11x_3)(22+  7x_1-5x_2+  8x_3)(19
 +  11x_1+  2x_2+  3x_3)\\
&(18+  11x_1+  2x_2-3x_3)(20+  8x_1+  7
 x_2+  5x_3),\\
\pp_{0,(1,2,1)}(x)=&(20-7x_1-5x_2-8x_3)(20-2x_1-3x_2-11x_3)(20-6x_1+
 -10x_2)\\
&(20+  6x_1-10x_2)(20+  8x_1-7x_2-5
 x_3)(21-2x_1+  3x_2-11x_3)\\
&(22-7x_1+  5x_2-8x_3)(22+  3
 x_1+  11x_2+  2x_3)(24-3x_1+  11x_2+  2x_3)\\
&(20+  7x_1+  5x_2+
 -8x_3)(21-11x_1-2x_2-3x_3)(19-8x_1-7x_2+  5x_3)\\
&(20+
 -11x_1-2x_2+  3x_3)(22+ 6x_2+  10x_3)(22-5x_1+  8x_2
 +  7x_3)\\
&(20-6x_2+  10x_3)(21+  5x_1-8x_2+  7x_3)(21
 +  10x_1+  6x_3)\\
&(19+  10x_1-6x_3)(22+  7x_1+  5
 x_2+  8x_3),\\
\pp_{0,(2,1,1)}(x)=
&(20-5x_1-8x_2-7x_3)(20-6x_2-10x_3)(20-3x_1-11x_2-2x_3)\\
&(20+3x_1-11x_2-2x_3)(20+  5x_1-8x_2-7x_3)(22+6x_2-10x_3)\\
&(23-5x_1+  8x_2-7x_3)(22+3x_1+11x_2-2x_3)(24-3x_1+11x_2-2x_3)\\
&(21+5x_1+8x_2-7x_3)(21-10x_1-6x_3)(18-7x_1-5x_2+8x_3)\\
&(19-10x_1+  6x_3)(20-2x_1+  3x_2+  11x_3)(20-7x_1+  5x_2+  8x_3)\\
&(19-2x_1-3x_2+  11x_3)(21+  8x_1-7x_2+  5x_3)(20+  11x_1-2x_2+  3x_3)\\
&(19+  11x_1-2x_2-3x_3)(22+  5x_1+  8x_2+  7x_3).
\end{split}
\end{equation*} 
}
Next we have to determine an $\overline{\epsilon}$ as defined in
\eqref{eq:epsilon}. To this end we have estimated all the needed distances
$\epsilon_k$ (cf.~\eqref{eq:epsilon_k}) by a rather
ad hoc method and found   
\begin{equation*}
  \epsilon_2> 1/10,\quad \epsilon_1> 3/50, \quad\epsilon_0> 3/100.
\end{equation*}
Hence we may set $\overline{\epsilon}=3/100$ and since
$\diam(P)\leq 4$ we may choose for the exponent $p$ of Lemma
\ref{lem:approximating} $p=332$. With these values we get from Lemma
\ref{lem:approximating} the 
following polynomial $\pp_{\overline{\epsilon}}(x)$  
{\tiny
\begin{equation*}
\begin{split}
 \pp_{\overline{\epsilon}}(x)=&\frac{1}{12}
      \left[\frac{6x_2+4x_3+1}{11}\right]^{664}
      +\frac{1}{12}
      \left[\frac{-6x_2+4x_3-1}{11}\right]^{664}
      +\frac{1}{12}
      \left[\frac{4x_1+6x_2+1}{11}\right]^{664} \\
      +&\frac{1}{12}
      \left[\frac{4x_1-6x_3+1}{9}\right]^{664}
      +\frac{1}{12}
      \left[\frac{6x_1+4x_2-1}{9}\right]^{664}
      +\frac{1}{12}
      \left[\frac{-6x_1+4x_2+1}{11}\right]^{664}\\
+&\frac{1}{12}
      \left[\frac{-6x_2-4x_3-1}{11}\right]^{664}
      +\frac{1}{12}
      \left[\frac{6x_2-4x_3+1}{11}\right]^{664}
      +\frac{1}{12}
      \left[\frac{-4x_1-6x_2-1}{11}\right]^{664} \\
      +&\frac{1}{12}
      \left[\frac{-4x_1+6x_3-1}{9}\right]^{664}
      +\frac{1}{12}
      \left[\frac{-6x_1-4x_2+1}{9}\right]^{664}
      +\frac{1}{12}
      \left[\frac{6x_1-4x_2-1}{11}\right]^{664}.
\end{split}
\end{equation*}
}

\section{Outlook}

Why should anyone care about the representation of polyhedra by
exponentially many polynomial inequalities, given that one knows
that quadratically many suffice? Our answer is that the latter result
is of pure existential nature, while we can {\em construct} such inequalities.
We admit that the representations we found do not form an achievement
of concrete practical value. That is why we did not state them
in an algorithmic fashion. We see our paper just as a small
step towards a development of real algebraic geometry
in a constructive direction. There are a number of possible
routes. We want to mention briefly what we are interested
in and what might be achievable.

It would be nice to have efficient (in a sense that can be made precise)
algorithms that provide, e.g., for polytopes $P$ given in the form of a
$\mathcal{V}$- or $\mathcal{H}$-representation, a
$\mathcal{P}$-representation $P=\newsymb(\pp_1,\dots,\pp_l)$ with a
number $l$ of polynomials that is polynomially bounded in the
dimension of the polytope. 
It may also be useful to be able to construct a small number of 
``simple'' polynomials $\pp_1,\dots,\pp_k$ 
such that $\newsymb(\pp_1,\dots,\pp_k)$ approximates $P$ well. Of
course, one can study similar problems concerning the representation
of arbitrary semi-algebraic sets. E.g., given a semi-algebraic set
$\mathcal{S}$, can it be represented by a system of polynomials with
total degree at most $k$, say? Can such a system be constructed
efficiently? How well can $\mathcal{S}$ be approximated by polynomials
of degree $k$? For polyhedra $P$ we know that there exists a
representation by polynomials of total degree $1$, but what can we say
about the minimum number of polynomials of degree $k$ representing $P$?   
 
To indicate  possible outcomes that may result from such a change of
representation, 
let us look at the very successful polyhedral approach
to combinatorial optimization. The basic idea here is to represent
combinatorial objects (such as the tours of a travelling salesman,
the independent sets of a matroid, or the stable sets in a graph)
as the vertices of a polytope. This way one arrives at an
(implicit) $\mathcal{V}$-representation of classes of polytopes such as
travelling salesman or stable set polytopes. If one can find complete or
tight partial representations of polytopes of this type by linear
equations and inequalities (i.e., $\mathcal{H}$-representations),
linear programming (LP) techniques can 
be employed to solve the associated combinatorial optimization
problem, see \cite{GLS:ellipsoid}. 

This approach provides a general machinery to establish the polynomial
time solvability of combinatorial problems theoretically. In
particular, it is often employed to identify easy special cases of
generally hard problems. One such example is the stable set
problem that is \NP-hard for general graphs but solvable in
polynomial time for perfect (and other classes of) graphs, see
\cite{GLS:ellipsoid}, Chapter 9.

The LP approach provides more. Even in the case where only
partial $\mathcal{H}$-representations of the polyhedra associated with
  combinatorial
problems are known, LP techniques (such as cutting planes and column
generation) have resulted in very successful exact or approximate
solution methods. One prime example for this methodology is the
travelling salesman problem, see \cite{applebixbyetal:travel} and the corresponding
web page at {\tt http://www.math.princeton.edu/tsp/}, which includes an annotated bibliography
with remarks about the historical development of this area.

Progress of the type  may also be possible via the
``$\mathcal{P}$-representation approach''. Let us discuss
this by means of the stable set problem.

Although complexity theory suggests that it is inconceivable that
one can find an explicit $\mathcal{P}$-representation of all members of the
class
of stable set polytopes, i.e., the convex hull of all incidence vectors
of stable sets, it might be possible to find, for every graph $G$,
a "small" number $l(G)$ of not too "ugly" polynomials such that
$\newsymb(\pp_1,\dots,\pp_{l(G)})$ approximates the stable
set polytope $\STAB(G)$ well, and such that
$\newsymb(\pp_1,\dots,\pp_{l(G)})$ equals $\STAB(G)$ for a special
class of graphs.

It is also conceivable that, for such particular systems
$\pp_1,\dots,\pp_{l(G)}$
of polynomials, special nonlinear programming algorithms can be
designed that solve optimization problems over
$\newsymb(\pp_1,\dots,\pp_{l(G)})$ efficient in practice or theory.

We do know, of course, that these indications of possible future
results are mere speculation. Visions of this type, however,
were the starting point of the results presented in this
paper. And we do hope that there will be progress in some of
the directions mentioned.

\vspace{0.5cm}
\noindent{\it Acknowledgment.} We would like to thank Monika Ludwig
and G\"unter M. Ziegler for helpful discussions, and
Michael Joswig for pointing out a mistake  in a previous version of
Corollary \ref{cor:polyhedra} and for giving suggestions for
improvements.  We are grateful to the  anonymous
referees for their helpful comments.

\providecommand{\bysame}{\leavevmode\hbox to3em{\hrulefill}\thinspace}

\end{document}